\input amstex

\documentstyle{amsppt}

\refstyle{A}

\nologo


\hoffset .25 true in
\voffset .2 true in

\hsize=6.1 true in
\vsize=8.2 true in

\define\Adot{\bold A^\bullet}

\define\dm{\operatorname{dim}}

\define\kr{\operatorname{ker}}

\define\cokr{\operatorname{coker}}

\define\image{\operatorname{im}}

\topmatter

\title Semi-simple Carrousels and the Monodromy \endtitle

\author David B. Massey \endauthor

\address{David B. Massey, Dept. of Mathematics, Northeastern University, Boston,
MA, 02115, USA} \endaddress

\email{dmassey\@neu.edu}\endemail


\abstract    Let $\Cal U$ be an open neighborhood of the origin in $\Bbb C^{n+1}$ and let $f:(\Cal U, \bold 0)\rightarrow(\Bbb C, 0)$ be complex analytic. Let $z_0$ be a generic linear form on $\Bbb C^{n+1}$. If the relative polar curve $\Gamma^1_{f, z_0}$ at the origin is irreducible and the intersection number $\big(\Gamma^1_{f, z_0}\cdot V(f))_\bold 0$ is prime, then there are severe restrictions on the possible degree $n$ cohomology of the Milnor fiber at the origin. We also obtain some interesting, weaker, results when $\big(\Gamma^1_{f, z_0}\cdot V(f))_\bold 0$ is not prime. 

\endabstract

\endtopmatter

\document

\noindent\S0. {\bf Introduction}  

\vskip .1in

In [{\bf L\^e2}] and [{\bf L\^e3}], L\^e introduces his carrousel as a tool for  analyzing the relative monodromy of the Milnor fiber of a function, $f$, modulo a hyperplane slice. In [{\bf T1}] and [{\bf T2}], Tib\u ar gives a careful presentation of L\^e's carrousel and uses it to obtain interesting results. Outside of the work of L\^e and Tib\u ar, the carrousel seems to be a largely unused device. This is due in part to the complicated nature of the carrousel description.

In this short paper, we look at some interesting special cases that occur and, in particular, look at the case where the relative polar curve, $\Gamma^1_{f, z_0}$, has a single component such that the intersection number $\big(\Gamma^1_{f, z_0}\cdot V(f)\big)_\bold 0$ is prime. In this case, we show, in Theorem 2.3, how L\^e's carrousel tells one a great deal about the middle-dimensional homology/cohomology groups of the Milnor fiber of $f$, regardless of the dimension of the critical locus.

\vskip .3in

\noindent\S1. {\bf L\^e's Playground}  

\vskip .1in

Let $\Cal U$ be an open neighborhood of the origin in $\Bbb C^{n+1}$ and let $f:(\Cal U, \bold 0)\rightarrow(\Bbb C, 0)$ be a complex analytic function which has a critical point at the origin.

Recall that a {\it good stratification for $f$} is a stratification $\Cal S$  of $V(f)$ which contains $V(f)-\Sigma f$, and such that, for all $S\in \Cal S$, the pair $(\Cal U-V(f), S)$ satisfies the $a_f$ condition. After a linear change of coordinates, we may assume that the first coordinate, $z_0$, is a prepolar form (or coordinate) for $f$ at $\bold 0$ (see [{\bf M1}]); this means that there exists a neighborhood, $\Cal W\subseteq\Cal U$, of $\bold 0$  such that, inside $\Cal W-\{\bold 0\}$, $V(z_0)$ transversely intersects all of the strata of a good stratification of $V(f)$ (we do {\bf not} need the condition of the frontier here -- we could simply use a good partition). Then, at the origin, the relative polar curve $\Gamma^1_{f, z_0}$ (see [{\bf M1}])  is purely one-dimensional (or empty), and $\Gamma^1_{f, z_0}$ properly intersects both $V(f)$ and $V(z_0)$ (again, see [{\bf M1}]). We always consider $\Gamma^1_{f, z_0}$ with its cycle structure (see [{\bf M1}]). We assume that $\Cal U$ is small enough so that every component of $\Gamma^1_{f, z_0}$ passes through the origin.

Let $D$ be a component of the cycle $\Gamma^1_{f, z_0}$ (with either its reduced structure or its cycle structure). We have the following well-known formula, originally due to Teissier, 
$$
\big(D\cdot V(f)\big)_\bold 0 = \big(D\cdot V(z_0)\big)_\bold 0 + \left(D\cdot V\left(\frac{\partial f}{\partial z_0}\right)\right)_\bold 0.
$$
As $\bold 0$ is a critical point of $f$, it follows that $\left(D\cdot V\left(\frac{\partial f}{\partial z_0}\right)\right)_\bold 0>0$, and so 
$$n_D:=\big(D\cdot V(f)\big)_\bold 0 > \big(D\cdot V(z_0)\big)_\bold 0 =: m_D.\tag{$\dagger$}$$

\vskip .2in

\noindent{\bf L\^e's Attaching Theorem}

\vskip .1in

Let $B_\epsilon$ (resp., $\Bbb D_\delta$)  denote a closed ball of radius $\epsilon$ (resp., $\delta$) centered at the origin in $\Bbb C^n$ (resp., $\Bbb C$). Assume that $0<\eta\ll\delta\ll\epsilon\ll 1$. Let $\xi\in\Bbb C$ be such that $0<|\xi|\leqslant\eta$. Then, $$F_f:=\big(\Bbb D_\delta\times B_\epsilon\big)\cap f^{-1}(\xi)$$ is (up to homotopy) the Milnor fiber of $f$ at $\bold 0$, and $F_{f_0}:= V(z_0)\cap F_{f, \bold 0}$ is the Milnor fiber of $f_0:= f_{|_{V(z_0)}}$ at the origin. The main theorem of [{\bf L\^e1}] (see, also, [{\bf M1}]) is:

\vskip .3in

\noindent{\bf Theorem 1.1}. (L\^e) {\it The Milnor fiber $F_{f}$ is obtained from $F_{f_0}$ by attaching ${\tau_{f, z_0}}:=\big(\Gamma^1_{f, z_0}\cdot V(f)\big)_\bold 0$ $n$-handles ($n$-cells, up to homotopy).

Thus,
$H^{k}(F_{f}, F_{f_0})=0$ if $k\neq n$, there is an isomorphism $${\omega_{f, z_0}}: H^{n}(F_{f}, F_{f_0})@>\cong>>\bigoplus_D\Bbb Z^{n_D}\ \cong\ \Bbb Z^{\tau_{f, z_0}},$$
where $D$ ranges over the (possibly non-reduced) components of $\Gamma^1_{f, z_0}$,
 and there is a map on reduced, integral cohomology $\widetilde H^{n-1}(F_{f_0})@>\partial_{f, z_0}>>H^{n}(F_{f}, F_{f_0})$ such that $\kr\partial_{f, z_0}\cong\widetilde H^{n-1}(F_{f})$ and $\cokr\partial_{f, z_0}\cong\widetilde H^{n}(F_{f})$.}

\vskip .2in

We refer to the above result as {\it L\^e's Attaching Theorem}.

\vskip .3in

\noindent{\it Remark 1.2}. By the naturality of the Milnor monodromy, the map $\partial_{f, z_0}$ commutes with the respect Milnor monodromies on $\widetilde H^{n-1}(F_{f_0})$ and $H^{n}(F_{f}, F_{f_0})$. In particular, the image of $\partial_{f, z_0}$, $\image\partial_{f, z_0}$, is a free Abelian submodule which is invariant under the monodromy.

In addition, the main theorem of A'Campo in [{\bf A'C}] tells us that the trace of the monodromy action on $H^{n}(F_{f}, F_{f_0})$ is $0$. Therefore, the trace of the monodromy action on $\image\partial_{f, z_0}$ is negative the trace of the monodromy action on the free part of $\widetilde H^{n}(F_{f})$.

\vskip .1in

For all $k$, we denote the rank of $\widetilde H^k(F_f)$ (i.e., the $k$-th reduced Betti number) by $\tilde b_k(f)$. Thus, the rank of $\image\partial_{f, z_0}$ is $e_{f, z_0}:=\tau_{f, z_0}-\tilde b_n(f)$. We denote the characteristic polynomials of the monodromy action on $\image\partial_{f, z_0}$, on $H^{n}(F_{f}, F_{f_0})$, and on the free part of $\widetilde H^{n}(F_{f})$ (or on $\widetilde H_{n}(F_{f})$) by ${\operatorname{char}}_{\image\partial_{f, z_0}}(\lambda)$, ${\operatorname{char}}_{{\operatorname{rel}}_{f, z_0}}(\lambda)$, and ${\operatorname{char}}^{n}_f(\lambda)$, respectively. Of course, we have the equality
$${\operatorname{char}}_{{\operatorname{rel}}_{f, z_0}}(\lambda)= {\operatorname{char}}_{\image\partial_{f, z_0}}(\lambda)\cdot {\operatorname{char}}^{n}_f(\lambda).$$

\vskip .2in

\noindent{\bf The Swing}

\vskip .1in

In [{\bf L-P}], L\^e and Perron use the ``swing'' to more carefully analyze the image of the attaching map $\partial_{f, z_0}$ above. They do this in the case where $\dm_\bold 0\Sigma f =1$. However, in [{\bf M1}], we showed that their argument works regardless of the dimension of the critical locus. 

What the swing shows is:

\vskip .3in

\noindent{\bf Theorem 1.3}. (L\^e and Perron) {\it The rank of the free Abelian module $\image\partial_{f, z_0}$ is at least $\gamma^1_{f, z_0}:=\big(\Gamma^1_{f, z_0}\cdot V(z_0)\big)_\bold 0$. Thus, the rank of $\widetilde H^{n}(F_{f})$ is at most 
$$\lambda^0_{f, z_0}:= {\tau_{f, z_0}}-\gamma^1_{f, z_0} = \left(\Gamma^1_{f, z_0}\cdot V\left(\frac{\partial f}{\partial z_0}\right)\right)_\bold 0.
$$

In fact, for each component $D$ of $\Gamma^1_{f, z_0}$, there is a submodule $E_D\subseteq\Bbb Z^{n_D}$ which is generated by $m_D$ of the basis elements of $\Bbb Z^{n_D}$ such that, if $\pi$ denotes the projection from $\bigoplus_D\Bbb Z^{n_D}$ onto $\bigoplus_D E_D\ \cong \Bbb Z^{\gamma^1_{f, z_0}}$, then $\pi\circ{\omega_{f, z_0}}\circ\partial_{f, z_0}$ is a surjection (where $\partial_{f, z_0}$ and ${\omega_{f, z_0}}$ are defined in Theorem 1.1).
}

\vskip .3in

The number $\lambda^0_{f, z_0}$ defined above is the {\it $0$-th L\^e number of $f$ (at the origin with respect to $z_0$)} (see [{\bf M1}]). 

\vskip .2in

\noindent{\bf L\^e's Monodromy Carrousel}

\vskip .1in

L\^e's Carrousel (see [{\bf L\^e2}] and [{\bf T2}]) gives a geometric description of the monodromy action on $H^{n}(F_{f}, F_{f_0})$. Let us briefly recall the set-up and some features of L\^e's Carrousel.

Let $\Theta$ denote the map $(z_0, f)$ from $\big(\Bbb D_\delta\times B_\epsilon\big)\cap f^{-1}(\Bbb D_\eta)$ onto $\Bbb D_\delta\times \Bbb D_\eta$, and use $(u, v)$ for coordinates on the codomain. Then, $C_{f, z_0}:= \Theta\big(\Gamma^1_{f, z_0}\big)$ is the {\it Cerf diagram} of $f$ with respect to $z_0$. Each component of $C_{f, z_0}$ is tangent to the $u$-axis at the origin (this follows from $(\dagger)$ above). The map $\Theta_{|_{\Gamma^1_{f, z_0}}}$ is finite, and we endow $C_{f, z_0}$ with a cycle structure via the proper push-forward. It follows that 
$$
{\tau_{f, z_0}}=\big(\Gamma^1_{f, z_0}\cdot V(f)\big)_\bold 0=\big(C_{f, z_0}\cdot V(v)\big)_\bold 0
$$
and
$$
\gamma^1_{f, z_0}=\big(\Gamma^1_{f, z_0}\cdot V(z_0)\big)_\bold 0=\big(C_{f, z_0}\cdot V(u)\big)_\bold 0.
$$

However, L\^e's Carrousel description requires that  $\Gamma^1_{f, z_0}$ be reduced and that $\Theta_{|_{\Gamma^1_{f, z_0}}}$ be one-to-one, i.e., that the cycle $C_{f, z_0}$ is reduced, and we have {\bf not} assumed that $z_0$ is generic enough to make this happen.

\vskip .3in

\noindent {\bf Definition 1.4}. The linear form $z_0$ is a {\it carrousel form (for $f$ at $\bold 0$)} if and only if $z_0$ is a prepolar form for $f$ at $\bold 0$ and, at the origin, the cycle $C_{f, z_0}$ is reduced.

\vskip .3in

{\bf We assume throughout the remainder of this section that $z_0$ is a carrousel form.}

\vskip .2in

Suppose, again, that $D$ is an irreducible component of $\Gamma^1_{f, z_0}$. Then, $C:=\Theta(D)$ is an irreducible component of $C_{f, z_0}$, and every component of $C_{f, z_0}$ is obtained in this manner. Moreover,
$$
m_C:= \big(C\cdot V(u)\big)_\bold 0 = \big(D\cdot V(z_0)\big)_\bold 0 = m_D ,
$$
and
$$
n_C:= \big(C\cdot V(v)\big)_\bold 0 = \big(D\cdot V(f)\big)_\bold 0 = n_D.
$$
Let $g_C:= \operatorname{gcd}(m_C, n_C)$, let $p_C:=m_C/g_C$, and let $q_C:= n_C/g_C$. The curve $C$ has a local parameterization of the form
$$v=t^{{}^{n_C}},\hskip .2in u= \alpha_{{}_C}t^{{}^{m_C}}+\text{ higher order terms }.$$
The {\it carrousel approximation of $C$} is the curve given by 
$$v=t^{{}^{q_C}},\hskip .2in u= \alpha_{{}_C}t^{{}^{p_C}},$$
i.e., the curve $\widehat C:=V(u^{{}^{q_C}}- \alpha_{{}_C}^{{}^{q_C}}v^{{}^{p_C}})$. We refer to $\beta_{{}_C}:=\alpha_{{}_C}^{{}^{q_C}}$ as the {\it carrousel coefficient of $C$}.

\vskip .3in

\noindent {\bf Definition 1.5}. We say that {\it carrousel of $f$ with respect to (the carrousel form) $z_0$ is semi-simple} provided that:

\vskip .1in

\noindent i)\hskip .14in for all components $C$ of $C_{f, z_0}$, $m_C$ and $n_C$ are relatively prime (and so, $(m_C, n_C)=(p_C, q_C)$);

\vskip .1in

\noindent ii)\hskip .1in distinct components of $C_{f, z_0}$ have distinct carrousel approximations, i.e., if $C_1\neq C_2$, then $(p_{C_1}, q_{C_1}, \beta_{C_1})\neq(p_{C_2}, q_{C_2}, \beta_{C_2})$.

\vskip .3in

 Before we state the next theorem, we need to give one more piece of terminology. We refer to the automorphism of $\Bbb Z^k$ given by $(a_1, a_2, \dots, a_k)\mapsto (-a_k, a_1, a_2, \dots, a_{k-1})$ as {\it cyclic anti-permutation}. The characteristic polynomial of cyclic anti-permutation is $\lambda^k+1$. Cyclic anti-permutation has $-1$ as an eigenvalue if and only if $k$ is odd and, in this case, the {\it anti-diagonal} $\widehat\Delta:=\Bbb Z(1, -1, 1, -1, \dots, -1, 1)$ is the eigenspace of $-1$. One shows easily that $\Bbb Z^k/\widehat\Delta\cong \Bbb Z^{k-1}$.
 
 \vskip .2in
 
 We use the terminology {\it semi-simple} because L\^e's carrousel study in [{\bf L\^e2}] immediately implies:

\vskip .3in

\noindent{\bf Theorem 1.6}. (L\^e) {\it If the carrousel of $f$ with respect to  $z_0$ is semi-simple, then the isomorphism of Theorem 1.1, $${\omega_{f, z_0}}: H^{n}(F_{f}, F_{f_0})@>\cong>>\bigoplus_D\Bbb Z^{n_D},$$
can be chosen so that each direct summand is invariant under the Milnor monodromy (i.e., the monodromy breaks up into blocks), and the action on each block is either cyclic permutation or cyclic anti-permutation. In particular,  the characteristic polynomial of the monodromy action on $H^{n}(F_{f}, F_{f_0})$ is 
$${\operatorname{char}}_{\operatorname{rel}}(\lambda)= \prod_D(\lambda^{n_D}\pm1).$$
}

\vskip .2in

\noindent{\it Proof}. One refers to the proofs in [{\bf L\^e2}]. 

Condition i) of being semi-simple implies, for a given component $C$ of $C_{f, z_0}$, that each carrousel disk contains at most one point of $C$; this implies that in one ``turn of the carrousel'' there is no interaction between different points on $C$.

Condition ii) of being semi-simple says that distinct components of $C_{f, z_0}$ have distinct carrousel approximations and, hence, the carrousel points of distinct components do not interact as the carrousel turns.

Now, the carrousel disks are permuted cyclically by the monodromy, and each carrousel disk centered at a point of a Cerf component $C$ contributes one copy of $\Bbb Z$ as a direct summand in $\Bbb Z^{n_C}=\Bbb Z^{n_d}$. However, after a carrousel disk returns to itself after $n_C=n_D$ iterations, the corresponding copy of $\Bbb Z$ may be mapped to itself by either plus or minus the identity. Hence, the induced map on cohomology is either cyclic permutation or anti-permutation, and the conclusion follows immediately.\qed

\vskip .5in

\noindent\S2. {\bf Prime Polar Curves}  

\vskip .1in

In this section, we continue with our notation from Section 1, and {\bf we continue to assume that $z_0$ is prepolar for $f$ at $\bold 0$, but  we no longer assume that $z_0$ is a carrousel form}.

\vskip .3in

\noindent{\bf Definition 2.1}. Let $D$ be a (possibly non-reduced) component  of the cycle $\Gamma^1_{f, z_0}$.

\vskip .1in

$D$ is called  {\it relatively prime} provided that $(D\cdot V(z_0))_{\bold 0}$ and $(D\cdot V(f))_{\bold 0}$ are relatively prime.

\vskip .1in

$D$ is called {\it unitary} provided that $(D\cdot V(z_0))_{\bold 0} = 1$, i.e., $D$ is reduced, $D$ is smooth, and $D$ is transversely intersected at $\bold 0$ by $V(z_0)$.

\vskip .1in

$D$ is called {\it prime of order $\frak p$} provided that $\frak p:=(D\cdot V(f))_{\bold 0}$ is a prime number.

\vskip .1in

Note that, as $(D\cdot V(z_0))_{\bold 0}<(D\cdot V(f))_{\bold 0}$, unitary and prime components are also relatively prime.

\vskip .1in

The cycle $\Gamma^1_{f, z_0}$ is itself said to be {\it relatively prime} (resp., {\it unitary},   resp., {\it prime of order $\frak p$}) provided that $\Gamma^1_{f, z_0}$ has one irreducible component and that component is relatively prime (resp., unitary, resp., prime of order $\frak p$). Note that if $\Gamma^1_{f, z_0}$ is unitary or prime, then it is relatively prime.

\vskip .3in

\noindent{\bf Proposition 2.2}. {\it Suppose that every component of $\Gamma^1_{f, z_0}$ is relatively prime. Then, $z_0$ is a carrousel form. In particular, $\Gamma^1_{f, z_0}$ is reduced.

\vskip .1in

If  $\Gamma^1_{f, z_0}$ is itself relatively prime, then the carrousel of $f$ with respect to $z_0$ is semi-simple.
}

\vskip .2in

\noindent{\it Proof}. Suppose that $D$ is an irreducible component of $\Gamma^1_{f, z_0}$, which has a possibly non-reduced cycle structure. Let $C$ be the proper push-forward of $D$ by $\Theta$, i.e., $C:=\Theta_*(D)$. Then,
$\big(C\cdot V(u)\big)_\bold 0 = \big(D\cdot V(z_0)\big)_\bold 0$,
and
$\big(C\cdot V(v)\big)_\bold 0 = \big(D\cdot V(f)\big)_\bold 0$. We must show that $C$ is reduced.

Suppose that, as cycles, $C= kC^\prime$, where $C^\prime$ is reduced. Then $k$ must divide both $\big(C\cdot V(u)\big)_\bold 0$ and $\big(C\cdot V(v)\big)_\bold 0$, which are relatively prime. Thus, $k=1$, and $z_0$ is a carrousel form.

\vskip .1in

The remaining claim follows immediately from the definition of a semi-simple carrousel.\qed

\vskip .5in

We wish to recall now the notion of the suspension of $f$ (see, for instance, [{\bf M1}]and [{\bf M3}]). Suppose that $f=f(\bold z)$. Suppose that $w$ is a variable disjoint from $\bold z$. Then, the function $f+w^2$ on $\Cal U\times\Bbb C$ is called the {\it suspension of $f$}. 

It is trivial to show that $\Sigma(f+w^2)=\Sigma f\times\{0\}$, that $\Gamma^1_{f+w^2, z_0}=\Gamma^1_{f, z_0}\times\{0\}$, and that if $z_0$ is prepolar for $f$ at $\bold 0$, then $z_0$ is prepolar for $f+w^2$ at $\bold 0$. See [{\bf M1}]. It follows easily that $\gamma^1_{f, z_0}= \gamma^1_{f+w^2, z_0}$ and $\tau_{f, z_0}= \tau_{f+w^2, z_0}$. Therefore, $\Gamma^1_{f, z_0}$ is  prime of order $\frak p$ if and only if $\Gamma^1_{f+w^2, z_0}$ is  prime of order $\frak p$ 

By the Sebastiani-Thom result (for references to this result, in many various cases, see [{\bf M1}] and [{\bf M3}]), for all $k$, $\widetilde H^{k+1}(F_{f+w^2})\cong \widetilde H^k(f_f)$ and, under this isomorphism, the Milnor monodromy action on $\widetilde H^{k+1}(F_{f+w^2})$ is negative the monodromy action on $\widetilde H^k(F_f)$. one then recovers (as we saw above) the isomorphism 
$$\widetilde H^{n+1}(F_{f+w^2}, F_{f_0+w^2})\cong \widetilde H^{n+1}(F_{f}, F_{f_0})
$$
and finds that, under this isomorphism, the Milnor monodromy action on $\widetilde H^{n+1}(F_{f+w^2}, F_{f_0+w^2})$ is negative the monodromy action on $\widetilde H^{n+1}(F_{f}, F_{f_0})$. Thus, we have the following relationships between characteristic polynomials of the Milnor monodromy actions:
$$
{\operatorname{char}}_{\image\partial_{f+w^2, z_0}}(\lambda)=(-1)^{e_{f, z_0}}\cdot{\operatorname{char}}_{\image\partial_{f, z_0}}(-\lambda)
$$

$${\operatorname{char}}_{{\operatorname{rel}}_{f+w^2, z_0}}(\lambda)=(-1)^{\tau_{f, z_0}}\cdot{\operatorname{char}}_{{\operatorname{rel}}_{f, z_0}}(-\lambda)
,$$
and
$${\operatorname{char}}_{f+w^2}^{n+1}(\lambda)=(-1)^{\tilde b_n(f)}\cdot{\operatorname{char}}_f^n(-\lambda).$$

\vskip .5in

Below, we state a result in terms of the homology of $F_f$, instead of cohomology. While, in general, we prefer to think in cohomological terms, discussions of the monodromy action on $\widetilde H^n(F_f)$ are more complicated by the possible presence of torsion. However, $\widetilde H_{n}(F_{f})$ is free Abelian and is thus isomorphic to the free part of $\widetilde H^n(F_f)$.

\vskip .3in

\noindent{\bf Theorem 2.3}. {\it Suppose that $\Gamma^1_{f, z_0}$ is  prime of order $\frak p$. Then, we are in one of the following non-overlapping cases:

\vskip .1in

\noindent Case 0: $\widetilde H_{n}(F_{f}) = 0$, $\operatorname{rank}\widetilde H_{n-1}(F_{f_0})\geqslant\frak p$, and $\dm_\bold 0\Sigma f\geqslant 1$;

\vskip .1in

\noindent Case 1: $\widetilde H_{n}(F_{f})\cong \Bbb Z$, and the monodromy action on $\widetilde H_{n}(F_{f})$ is either \hskip .05in {\rm a)} the identity or \hskip .05in {\rm b)} negative the identity;

\vskip .1in

\noindent Case 2: $\frak p\neq 2$, $\widetilde H_{n}(F_{f})\cong \Bbb Z^{\frak p-1}$, $\widetilde H_{n-1}(F_{f})$ is free Abelian,  $\Gamma^1_{f, z_0}$ is unitary, and the characteristic polynomial of the monodromy action on $\widetilde H_{n}(F_{f})$ is either  \hskip .05in {\rm a)} $(\lambda^{\frak p}-1)/(\lambda-1)$ or  \hskip .05in {\rm b)} $(\lambda^{\frak p}+1)/(\lambda+1)$.

\vskip .1in

Moreover,  if $\frak p=2$ and $\widetilde H_{n}(F_{f})\cong \Bbb Z$, then $\widetilde H_{n-1}(F_{f})$ is free Abelian,  and $\Gamma^1_{f, z_0}$ is unitary.

\vskip .1in

In addition, suspending $f$ (and using the ``same'' coordinate $z_0$) leaves one in the same case, but interchanges the subcases {\rm a)} and {\rm b)} in Cases 1 and 2. }

\vskip .2in

\noindent{\it Proof}. By Proposition 2.2, $z_0$ is a carrousel form and the carrousel of $f$ with respect to $z_0$ is semi-simple. Therefore, by Theorem 1.6, either ${\operatorname{char}}_{{\operatorname{rel}}_{f, z_0}}(\lambda)= \lambda^{\frak p}-1$ or ${\operatorname{char}}_{{\operatorname{rel}}_{f, z_0}}(\lambda)= \lambda^{\frak p}+1$.

By 1.3, $\image\partial_{f, z_0}$ is non-zero. Therefore, Remark 1.2 implies that $\widetilde H_{n}(F_{f}) = 0$, or that ${\operatorname{char}}_f^n(\lambda)$ is $\lambda-1$, $\lambda+1$, $(\lambda^{\frak p}-1)/(\lambda-1)$, or $(\lambda^{\frak p}+1)/(\lambda+1)$, where this last characteristic polynomial cannot occur if $\frak p=2$.

\vskip .1in

Case 0: Suppose that  $\widetilde H_{n}(F_{f}) = 0$; this is equivalent to $\operatorname{rank}(\image\partial_{f, z_0})=\frak p$. This certainly implies that $\operatorname{rank}\widetilde H_{n-1}(F_{f_0})\geqslant\frak p$. We claim that it follows that $f$ cannot have an isolated critical at the origin. Suppose to the contrary that $\dm_\bold 0\Sigma f=0$. Then, by the formula of L\^e and Greuel, $\mu_\bold 0(f)+ \mu_\bold 0(f_0) = \frak p$, where $\mu$ denotes the Milnor number. As $\mu_\bold 0(f)>0$, $\mu_\bold 0(f_0) = \operatorname{rank}\widetilde H_{n-1}(F_{f_0})<\frak p$, and we are finished.

\vskip .1in

Case 1: ${\operatorname{char}}_f^n(\lambda)=\lambda\pm 1$. The claims follow immediately.

\vskip .1in

Case 2: ${\operatorname{char}}_f^n(\lambda)$ is  $(\lambda^{\frak p}-1)/(\lambda-1)$ or $(\lambda^{\frak p}+1)/(\lambda+1)$. Then, ${\operatorname{char}}_{\image\partial_{f, z_0}}(\lambda)$ is $\lambda\pm 1$; thus, under the isomorphism ${\omega_{f, z_0}}$, $\image\partial_{f, z_0}$ is contained in the diagonal or anti-diagonal of $\Bbb Z^{\frak p}$. By the first statement of Theorem 1.3, it follows that $\left(\Gamma^1_{f, z_0}\cdot V(z_0)\right)_\bold 0=1$  and so $\Gamma^1_{f, z_0}$ is unitary. Now, the last statement of Theorem 1.3 implies that $\image\partial_{f, z_0}$ must be the entire diagonal or anti-diagonal. It follows that $\operatorname{coker}\partial_{f, z_0}\cong \widetilde H^n(F_f)$ is free Abelian, which is equivalent to $\widetilde H_{n-1}(F_f)$ being free Abelian.

\vskip .1in

The suspension claim is immediate from the properties discussed prior to the theorem.\qed

\vskip .3in

\noindent{\it Example 2.4}. We will show here that all of the cases of Corollary 2.3 can occur. 

\vskip .1in

Note that, if $\dm_\bold 0\Sigma f\leqslant 1$, then $z_0$ is prepolar if and only if $\dm_\bold 0\Sigma f_0\leqslant 0$.

\vskip .1in

First, consider $f=z_0^2+ z_1^2+\dots+z_n^2$. Then, we know that $\widetilde H_{n-1}(F_{f_0})\cong\Bbb Z$ and $\widetilde H_{n}(F_{f})\cong\Bbb Z$. By A'Campo's main theorem in [{\bf A'C}], the trace of the monodromy action on $\widetilde H_{n}(F_{f})$ is $(-1)^{n+1}$.

Now, as a cycle,
$$\Gamma^1_{f, z_0}= V\left(\frac{\partial f}{\partial z_1}, \dots, \frac{\partial f}{\partial z_n}\right) = V(z_1, \dots, z_n).$$
Therefore, $\Gamma^1_{f, z_0}$ has a single component and ${\tau_{f, z_0}} = 2$ is prime, and so we can apply Theorem 2.3. By looking at the trace, we conclude that we are in Case 1a if $n$ is odd, and in Case 1b if $n$ is even.

\vskip .2in

Now, we will give examples of Case 2. Suppose that $\dm_\bold 0\Sigma f=0$, $\dm_\bold 0\Sigma f_0=0$, and that $\Gamma^1_{f, z_0}$ is  prime of order $\frak p\geqslant 3$. Then, $\widetilde H_n(F_f)\cong\Bbb Z^{{}^{\mu_{\bold 0}(f)}}$. Therefore, if $\mu_{\bold 0}(f)\geqslant 2$, then we must be in Case 2, and the trace distinguishes subcases a) and b). 

To give a specific example, let $f=y^2-x^3$, where we use $x$ in place of $z_0$. Then, $\mu_\bold 0(f)=2$, $\Gamma^1_{f, x} = V(y)$, and $\tau_{f, x}=3=\frak p$. By A'Campo's Lefschetz number result, we must be in Case 2b. By suspending, we find that $g:= w^2+y^2-x^3$ (again, using $z_0=x$) would be an example of Case 2a.

\vskip .2in

The example that we use for Case 0 was first shown to us by Dirk Siersma. Consider $f=(x^2+y^2-z^2)(y-z)$. We use the coordinate $z$ for $z_0$. The critical locus of $f$ is the line $V(x, y-z)$. As cycles, we find
$$
V\left(\frac{\partial f}{\partial x}, \frac{\partial f}{\partial y}\right) = V(2x(y-z), \ 2y(y-z)+x^2+y^2-z^2) = 3V(x, y-z)+ V(x, 3y+z).
$$
Therefore, $\Gamma^1_{f,z} = V(x, 3y+z)$ and ${\tau_{f, z_0}} = 3$. One also finds that $\mu_\bold 0(f_0)=4$. Thus, it is at least possible that we are in Case 1, but we must show this.

After an analytic coordinate change, $f=(x^2+st)t$. As $f$ is homogeneous, $F_f$ is diffeomorphic to $f^{-1}(1)$. Now we observe that $f^{-1}(1)$ is the set of points where $t\neq0$ and $s=(1-tx^2)/t^2$. Thus, $F_f$ is diffeomorphic to $\Bbb C\times\Bbb C^*$, and so is homotopy-equivalent to $S^1$. It follows that $\widetilde H_2(F_f)=0$.

\vskip .5in

We would like to show that Case 2 of Theorem 2.3 rarely occurs. For this, we will need the result below.

\vskip .3in

\noindent{\bf Proposition 2.5}. {\it Suppose  that the rank of $\widetilde H_n(F_f)$ equals $\lambda^0_{f, z_0}$.

\vskip .1in

Then,  the trace of the monodromy action on $\widetilde H_n(F_f)$ is $$(-1)^{n+1}\big(1-\chi(\Bbb L_{\Sigma f, z_0})\big),$$
where $\chi$ denotes the Euler characteristic and $\Bbb L_{\Sigma f, z_0}$ is the ``complex link of $\Sigma f$ at the origin with respect to $z_0$'', i.e., 
$$\Bbb L_{\Sigma f, z_0}:={\overset\circ\to B}_\epsilon\cap \Sigma f\cap V(z_0-\delta),$$ where $0\ll |\delta|\ll\epsilon\ll 1$ and ${\overset\circ\to B}_\epsilon$ is an open ball of radius $\epsilon$ centered at $\bold 0$. Thus, 

In particular, if $\Sigma f$ itself is smooth and transversely intersected by $V(z_0)$ at the origin, then the trace of the monodromy action on $\widetilde H_n(F_f)$ is $0$.
}

\vskip .2in

\noindent{\it Proof}. Recall from Remark 1.2 that the trace of the monodromy action on $\image\partial_{f, z_0}$ is negative the trace of the monodromy action on $\widetilde H_n(F_f)$.

Suppose that $\operatorname{rank}\widetilde H_n(F_f) = \lambda^0_{f, z_0}$.  In the case where $\dm_\bold 0\Sigma f \leqslant 1$,  the analysis of the ``nexus diagram'' in Application 2 of [{\bf M2}] tells us  that the trace of the monodromy action on $\image\partial_{f, z_0}$ is $$(-1)^{n+1}\big((|\Sigma f|\cdot V(z_0))_\bold 0-1\big)=(-1)^n\big(1-\chi(\Bbb L_{\Sigma f, z_0})\big).$$

As we commented at the end of [{\bf M2}], when the dimension of $\Sigma f$ is arbitrary, the nexus diagram still exists in the Abelian category of perverse sheaves, and the exact proof that we used when $\dm_\bold 0\Sigma f \leqslant 1$ tells us that there is an equality of Lefschetz numbers of the respective monodromy actions at the origin given by
$$
\Cal L_\bold 0\{\image\partial_{f, z_0}\}= \Cal L_\bold 0\{\phi_{f_0}[-1]\Bbb Z^\bullet_{V(z_0)}[n]\}-\Cal L_\bold 0\{\psi_{\hat z_0}[-1]\phi_f[-1]\Bbb Z^\bullet_{\Cal U}[n+1]\},
$$
where $\Cal L_\bold 0\{\Adot\}$ denotes the Lefschetz number at the origin of the Milnor monodromy action on the complex $\Adot$, and $\hat z_0$ is the restriction of $z_0$ to $V(f)$.
Now, $\image\partial_{f, z_0}$ is a sub-perverse sheaf of a perverse sheaf which is supported on a point; hence, $\Cal L_\bold 0\{\image\partial_{f, z_0}\}$ is simply the trace of the monodromy action on $\image\partial_{f, z_0}$. In addition, as we are assuming that $f$ has a critical point at the origin, A'Campo's result in [{\bf A'C}] implies that $\Cal L_\bold 0\{\phi_{f_0}[-1]\Bbb Z^\bullet_{V(z_0)}[n]\}=(-1)^n$. It remains for us to show that 
$$
\Cal L_\bold 0\{\psi_{\hat z_0}[-1]\phi_f[-1]\Bbb Z^\bullet_{\Cal U}[n+1]\} = (-1)^n\chi(\Bbb L_{\Sigma f, z_0}).\tag{$\dagger$}
$$

\vskip .1in

Consider the fundamental short exact sequence of perverse sheaves:
$$
0\rightarrow \Bbb Z_{V(f)}[n]\rightarrow \psi_f[-1]\Bbb Z^\bullet_{\Cal U}[n+1]\rightarrow \phi_f[-1]\Bbb Z^\bullet_{\Cal U}[n+1]\rightarrow 0.
$$

Let $\check z_0$ be the restriction of $z_0$ to $\Sigma f$. 
If we restrict this sequence to $\Sigma f$, then apply $\psi_{\check z_0}[-1]$,  and use that locally $\Sigma f\subseteq V(f)$, we obtain a distinguished triangle
$$
\psi_{\check z_0}[-1]\big(\Bbb Z_{\Sigma f}[n]\big)\rightarrow \psi_{\check z_0}[-1]\big(\left(\psi_f[-1]\Bbb Z^\bullet_{\Cal U}[n+1]\right)_{|_{\Sigma f}}\big)\rightarrow \psi_{\check z_0}[-1]\big(\left(\phi_f[-1]\Bbb Z^\bullet_{\Cal U}[n+1]\right)_{|_{\Sigma f}}\big)@>[1]>>,
$$
on which the monodromy acts compatibly. Using A'Campo's result again, we obtain that $$\Cal L_\bold 0\big(\psi_{\check z_0}[-1]\big(\left(\psi_f[-1]\Bbb Z^\bullet_{\Cal U}[n+1]\right)_{|_{\Sigma f}}\big)\big) = 0.$$
Thus, by additivity, we obtain that
$$
\Cal L_\bold 0\big(\psi_{\check z_0}[-1]\big(\left(\phi_f[-1]\Bbb Z^\bullet_{\Cal U}[n+1]\right)_{|_{\Sigma f}}\big)\big) = -\Cal L_\bold 0\big(\psi_{\check z_0}[-1]\big(\Bbb Z_{\Sigma f}[n]\big)\big).
$$
As the support of $\phi_f[-1]\Bbb Z^\bullet_{\Cal U}[n+1]$ already lies in $\Sigma f$, we obtain $(\dagger)$.
\qed

\vskip .3in

\noindent{\bf Corollary 2.6}. {\it Suppose that $\Gamma^1_{f, z_0}$ is  prime of order $\frak p$, and $\chi(\Bbb L_{\Sigma f, z_0})$ does not equal $0$ or $2$. Then, Case 2 of Theorem 2.3 does not occur, nor does Case 1 if $\frak p=2$.

In particular, if $\Gamma^1_{f, z_0}$ is  prime and $\Sigma f$ is itself smooth and transversely intersected by $V(z_0)$ at $\bold 0$, then the rank of $\widetilde H_n(F_f)$ is $0$ or $1$.
}

\vskip .2in

\noindent{\it Proof}. In Case 2 of Theorem 2.3, or in Case 1 if $\frak p=2$, the rank of $\widetilde H_n(F_f)$ equals $\tau_{f, z_0}-1=\lambda^0_{f, z_0}$, while the trace of the monodromy is $\pm 1$. The Corollary follows at once from Proposition 2.5.
\qed

\vskip .3in

In Example 2.4, we gave an example of a hypersurface with a line singularity which is a Case 0 example of Theorem 2.3. We also gave Case 1 examples which had isolated singularities. Corollary 2.6 tells us that we cannot produce a Case 2 example with a line singularity. Below, we give an example of a hypersurface with a line singularity which is Case 1.

\vskip .3in

\noindent{\it Example 2.7}. Consider the classic presentation of the Whitney umbrella as a family of nodes degenerating to a cusp: $f=y^2-x^3-tx^2$, where we use $t$ for our prepolar coordinate. Then, $\mu_\bold 0(f) =2$, $\Gamma^1_{f, t}= V(y, 3x+2t)$, $\tau_{f, t}= 3=\frak p$, and $\lambda^0_{f, t}=2$. Thus, up to isomorphism, $\partial_{f, t}$ is a map from $\Bbb Z^2$ to $\Bbb Z^3$. Therefore, $\operatorname{rank}\widetilde H_2(F_f)\geqslant 1$. However, as our critical locus is a line, we must be in Case b) of Theorem 2.5, and so $\operatorname{rank}\widetilde H_2(F_f)<2$. We conclude the well-known: $\widetilde H_2(F_f)\cong \Bbb Z$, i.e., this is an example of Case 1 of Theorem 2.3. 

\vskip .5in

\noindent\S3. {\bf More Complicated Examples}  

\vskip .1in

\noindent{\it Example 3.1}. Consider the family of examples $g(t,x,y):= y^2-x^a-t^cx^b$, where $a, b, c \geqslant 2$, and $a$ and $b$ are relatively prime. If $a\leqslant b$, then $g=y^2-x^a(1-t^cx^{b-a})$, which after an analytic change of coordinates at the origin becomes $y^2-x^a$; this is simply a cross-product of an isolated hypersurface singularity. So, assume that $a>b\geqslant 2$. We also assume that $a-b$ and $c$ are relatively prime. Note that this example subsumes Example 2.7.

\vskip .1in

One easily shows that $\Sigma g= V(x,y)$, and $g_0:=g_{|_{V(t)}}$ has an isolated critical point at the origin. Hence, $t$ is a prepolar coordinate for $g$.

Now, the Milnor number of $g_0$ at the origin is $a-1$ and, hence, the reduced cohomology of $F_{g_0}$ is $0$ in degree $0$ and is isomorphic to $\Bbb Z^{a-1}$ in degree $1$. We would like to know, ${\operatorname{char}}^{1}_{g_0}(\lambda)$, the characteristic polynomial of the monodromy action on $\widetilde H^1(F_{g_0})$.

The function $g_0$ is the suspension of the function $-x^a$ on $\Bbb C$. The Milnor fiber of $-x^a$ is $a$ points, which are permuted cyclically by the Milnor monodromy. Thus, the characteristic polynomial of the monodromy on the {\bf reduced} cohomology $\widetilde H^0(F_{-x^a})$ is $(\lambda^a-1)/(\lambda-1)$ and so, $${\operatorname{char}}^{1}_{g_0}(\lambda) = (\lambda^a-(-1)^{a})/(\lambda+1).$$

Now select a small $t_0\neq 0$. In the main theorem of [{\bf M2}], we proved that, if $\widetilde H^1(F_g)\neq 0$, then ${\operatorname{char}}^1_g(\lambda)$ not only divides  ${\operatorname{char}}^{1}_{g_0}(\lambda)$, but also divides ${\operatorname{char}}^1_{g_{t_0}}(\lambda)$, where $g_{t_0}$ denotes $g_{|_{V(t-t_0)}}$. 

Now, after an analytic change of coordinates at the origin (in $V(t-t_0)$), $g_{t_0}= y^2-x^b(x^{a-b}-t_0^c)$ becomes $y^2-x^b$. As this is the suspension of $-x^b$, we may use an analysis like that above to conclude that ${\operatorname{char}}^1_{g_{t_0}}(\lambda)$ equals $(\lambda^b-(-1)^{b})/(\lambda+1)$.

As $a$ and $b$ are relatively prime, we conclude that ${\operatorname{char}}^{1}_{g_0}(\lambda)$ and ${\operatorname{char}}^1_{g_{t_0}}(\lambda)$ have no common divisors. Therefore, we conclude that $\widetilde H^1(F_g)=0$.

\vskip .1in

As $V\left(\frac{\partial g}{\partial x}, \frac{\partial g}{\partial y}\right) = V(ax^{a-1}+bt^cx^{b-1}, y)$, we find that $\Gamma^1_{g, t} = V(ax^{a-b}+bt^c, y)$. Hence, $\gamma^1_{g, t} = \big(\Gamma^1_{g, t}\cdot V(t)\big)_\bold 0 =a-b$, and $$\lambda^0_{g,t} =  \Big(\Gamma^1_{g, t}\cdot V\left(\frac{\partial g}{\partial t}\right)\Big)_\bold 0 =  \big(\Gamma^1_{g, t}\cdot V(t^{c-1}x^b)\big)_\bold 0 = (c-1)(a-b)+bc = ac-(a-b).$$
Thus, $\tau_{g, t} = \gamma^1_{g, t}+\lambda^0_{g, t} = ac$. As $a$ and $b$ are relatively prime, so are $a-b$ and $a$. Therefore, as $a-b$ and $c$ are also relatively prime, we find that $\gamma^1_{g,t}$ and $\tau_{g,t}$ are relatively prime. In addition, since $a-b$ and $c$ are relatively prime, the polar curve $\Gamma^1_{g, t} = V(ax^{a-b}+bt^c, y)$ has a single irreducible component.

We conclude that $\Gamma^1_{g, t}$ is relatively prime and, hence, the carrousel of $g$ with respect to $t$ is semi-simple. Thus, ${\operatorname{char}}_{{\operatorname{rel}}_{g, t}}(\lambda) = \lambda^{ac}\pm 1$.

So, we have the map $\partial_{g,t}:\widetilde H^1(F_{g_0})\rightarrow H^2(F_{g}, F_{g_0})$, where $\widetilde H^1(F_{g_0})\cong\Bbb Z^{a-1}$, $H^2(F_{g}, F_{g_0})\cong \Bbb Z^{ac}$, ${\operatorname{char}}^{1}_{g_0}(\lambda) = (\lambda^a-(-1)^{a})/(\lambda+1)$, ${\operatorname{char}}_{{\operatorname{rel}}_{g, t}}(\lambda) = \lambda^{ac}\pm 1$, and $\widetilde H^1(F_g)=0$. It follows that $\widetilde H_2(F_g)\cong \Bbb Z^{ac-a+1}$ and that 
$${\operatorname{char}}^{2}_{g}(\lambda)= (\lambda^{ac}\pm 1)(\lambda+1)/(\lambda^a-(-1)^{a}).$$

\vskip .2in

Note that, if $b=2$, this example is an {\it isolated line singularity}, as studied by Siersma in [{\bf S}]. In this case, Siersma's work tells us a bit more: it says that $F_g$ has the homotopy-type of a bouquet of $(ac-a+1)$ $2$-spheres.

\vskip .3in

\noindent{\it Example 3.2}. In this example, we will look at  $f(s, t, x, y)=y^2-x^4+(s^3-t^2)x^3$. One easily checks that $\Sigma f = V(x, y)$, and so $f$ has a $2$-dimensional critical locus. Note also that $f_0:=f_{|_{V(s)}}$ is a function of the form of $g$ from Example 3.1, with $a=4$, $b=3$, and $c=2$.

\vskip .1in

We wish to see what our results can tell us about the cohomology and the monodromy of the Milnor fiber.

\vskip .1in

Our first problem is to verify that $s$ is a prepolar coordinate for $f$ at the origin. This means that we must first produce a good stratification. For this, we use the L\^e cycles and numbers, and apply  Corollary 6.6. and Remark 6.7 of [{\bf M1}]. We fix the coordinate system $(s, t, x, y)$ and will suppress any further reference to the coordinates.

We proceed with the calculation of the polar and L\^e cycles (see [{\bf M1}]):

\vskip .1in

$$
\Gamma^3_f= V\left(\frac{\partial f}{\partial y}\right) = V(y);
$$

$$
\Gamma^3_f\cdot V\left(\frac{\partial f}{\partial x}\right) = V(y)\cdot V\big(-4x^3+3(s^3-t^2)x^2\big) =
$$

$$
V(y)\cdot \Big(V\big(-4x+3(s^3-t^2)\big)\ +\ 2V(x)\Big) = V(-4x+3(s^3-t^2), y) + 2V(x, y)= 
$$

$$
\Gamma^2_f + \Lambda^2_f;
$$

$$
\Gamma^2_f\cdot V\left(\frac{\partial f}{\partial t}\right) = V(-4x+3(s^3-t^2), y)\cdot V(-2tx^3)=
$$

$$
V(-4x+3(s^3-t^2), y)\cdot \big(V(t)+3V(x)\big)=
$$

$$
V(-4x+3s^3, y, t) + 3V(s^3-t^2, x, y) = \Gamma^1_f + \Lambda^1_f;
$$

and, finally,

$$
\Gamma^1_f\cdot V\left(\frac{\partial f}{\partial s}\right) = V(-4x+3s^3, y, t)\cdot V(3s^2x^3) = 2[\bold 0]+3\cdot 3[\bold 0] = 11[\bold 0] = \Lambda^0_f.
$$

Thus,  we have $\Lambda^2_f = 2V(x, y)$, $\Lambda^1_f = 3V(s^3-t^2, x, y)$, and $\Lambda^0_f = 11[\bold 0]$. 

\vskip .1in

One easily calculates the L\^e numbers at a point $\bold p:=(s_0, t_0, x_0, y_0)$ near $\bold 0$:

\vskip .1in

\noindent $\lambda^0_f(\bold p)$ equals $11$ at the origin, and equals $0$ elsewhere;

\vskip .1in

\noindent $\lambda^1_f(\bold p) = \big(3V(s^3-t^2, x, y)\cdot V(s-s_0)\big)_{\bold p}$ equals $3\cdot 2=6$ at the origin, $3$ at other points of $V(s^3-t^2,x,y)$, and equals $0$ elsewhere;

\vskip .1in

\noindent $\lambda^2_f(\bold p)$ equals $2$ at all points of $V(x,y)$, and equals $0$ elsewhere.

\vskip .1in

Therefore, Corollary 6.6. of [{\bf M1}] tells us that $V(f)$ has a good stratification at the origin: $$\{V(f)-V(x,y),\ V(x,y)-V(s^3-t^2,x,y),\ V(s^3-t^2,x, y)-\{\bold 0\}, \{\bold 0\}\}.$$
Now,  $f_0= y^2-x^4-t^2x^3$ has a critical locus consisting of just the $t$-axis, i.e., $V(x,y)$ inside $V(s)$. We see then that $V(s)$ transversely intersects all of the good strata, except $\{\bold 0\}$, in a neighborhood of the origin, i.e., $s$ is prepolar for $f$ at $\bold 0$.

\vskip .1in

We continue to calculate (and continue to suppress the coordinates in the notation):

$$
\gamma^1_f = \big(\Gamma^1_f\cdot V(s)\big)_\bold 0 = \big(V(-4x+3s^3, y, t)\cdot V(s)\big)_\bold 0 = 1;
$$
and
$$
\tau_f = \gamma^1_f+\lambda^0_f = 1+11 = 12.
$$
As $\gamma^1_f=1$, $\Gamma^1_f$ is unitary.  Proposition 2.2 tells us that $s$ is a carrousel coordinate and that the carrousel of $f$ with respect to $s$ is semi-simple.

Putting all of the above work together, including our result in Example 3.1, we find that the map $\partial_{f, s}:\widetilde H^2(F_{f_0})\rightarrow H^3(F_f, F_{f_0})$ is a map from a $\Bbb Z$-module of rank $5$ to a copy of $\Bbb Z^12$, and the respective characteristic polynomials of the monodromy, acting on the free parts, are ${\operatorname{char}}^{2}_{f}(\lambda)= (\lambda^{8}\pm 1)(\lambda+1)/(\lambda^4-1)$ and ${\operatorname{char}}_{{\operatorname{rel}}_{f, s}}(\lambda)= \lambda^{12}\pm 1$.

Thus, in ${\operatorname{char}}^{2}_{f}(\lambda)$, we must choose the minus sign, and so ${\operatorname{char}}^{2}_{f_0}(\lambda)= (\lambda+1)(\lambda^4+1)$. By Proposition 2.5, the rank of $\image(\partial_{f, s})$ cannot be $1$. Therefore, we must have one of two cases:

\vskip .1in

\noindent i) $\widetilde H^0(F_f) = 0$, $\widetilde H^1(F_f) = 0$, $\operatorname{rank}\widetilde H^2(F_f) = 1$, $\operatorname{rank}\widetilde H^3(F_f) = 12-4=8$, ${\operatorname{char}}^{2}_{f}(\lambda)= \lambda+1$, and ${\operatorname{char}}^{3}_{f}(\lambda)= (\lambda^{12}+1)/(\lambda^4+1)$;

\vskip .1in

\noindent or

\vskip .1in

\noindent ii) $\widetilde H^0(F_f) = 0$, $\widetilde H^1(F_f) = 0$, $\operatorname{rank}\widetilde H^2(F_f) = 0$, $\operatorname{rank}\widetilde H^3(F_f) = 12-5=7$, and $${\operatorname{char}}^{3}_{f}(\lambda)= (\lambda^{12}+1)/[(\lambda+1)(\lambda^4+1)].$$

\vskip .1in

We do not, in fact, know which of these cases is the correct one.

\vskip .5in

\noindent\S4. {\bf Concluding Remarks}  

\vskip .1in

The main point of this paper is that the single number $\tau_{f, z_0}=\Big(\Gamma^1_{f, z_0}\cdot V(f)\Big)_\bold 0$ can tell one a great deal about $\widetilde H_n(F_f)$, at least when $\tau_{f, z_0}$ is prime.  

However, the calculation of $\tau_{f, z_0}$ is not a simple algebra exercise for complicated $f$.  If one wants to use Theorem 2.3, 
one must first ``calculate'' the polar curve, see that it has only one component, and then see that $\tau_{f, z_0}$ is prime. Moreover, as we saw in Example 3.2, if $\dm_\bold 0\Sigma f\geqslant 2$, then it is nontrivial to verify that $z_0$ is a prepolar coordinate.

Nonetheless, the case where $\Gamma^1_{f, z_0}$ is prime occurs in enough interesting examples that we find Theorem 2.3 to be interesting.

\vskip .3in

Together with L\^e D\~ung Tr\'ang, we believe that we can prove a generalization of Theorem 2.5. That generalization says that if $f$ has a smooth $1$-dimensional critical locus, and $z_0$ is a carrousel form, then either $\Gamma^1_{f, z_0}=\emptyset$ or $\operatorname{rank}\widetilde H_n(F_f)<\lambda^0_{f, z_0}$. This result requires a more detailed study of the carrousel and the swing.

\enddocument